\documentclass[12pt]{amsart}
\usepackage{amssymb}
\usepackage{txfonts}
\usepackage{amsmath}
\usepackage{amsfonts}
\usepackage{bbm}
\usepackage{mathrsfs}
\usepackage{wasysym}
\usepackage{dsfont}
\usepackage{hyperref}

\newtheorem{theorem}{Theorem}[section]
\newtheorem{lemma}[theorem]{Lemma}
\newtheorem{prop}[theorem]{Proposition}
\newtheorem{coro}[theorem]{Corollary}
\theoremstyle{definition}

\newtheorem{example}[theorem]{Example}

\theoremstyle{remark}
\newtheorem{remark}[theorem]{Remark}

\numberwithin{equation}{section} \hoffset=-.75in \textwidth=6in
\begin{document}

\title{subRiemannian geodesics of Carnot Groups of Step 3}

\author{Kanghai Tan}
\address{Department of Applied Mathematics, Nanjing University of Science \& Technology, Nanjing 210094, P.R. China}
\email{khtan@mail.njust.edu.cn}
\thanks{The first author is supported by NSF of China (No.10801073) and a grant from China Scholarship Council for study abroad.}

\author{Xiaoping Yang}
\address{School of Science, Nanjing University of Science \& Technology, Nanjing 210094, P.R. China}
\email{yangxp@mail.njust.edu.cn}


\date{}



\begin{abstract}
In Carnot groups of step$\leq3$, all subriemannian geodesics are
proved to be normal. The proof is based on  a reduction argument and
the Goh condition for minimality of singular curves. The Goh
condition is deduced from a reformulation and a calculus of the
end-point mapping which boils down to the graded structures of
Carnot groups.
\end{abstract}

\maketitle


\section{Introduction}
This paper is inspired by the smoothness problem of subriemannian
geodesics, one of the fundamental problems in subriemannian
geometry. We first study the case of Carnot group with step$\leq3$.
In this case we proved that all subriemannian geodesics are normal
and thus smooth.

To state the subriemannian geodesic problem, we first recall some
basic facts on subriemannian geometry. We refer to the book
\cite{Mo1} for detail. A subriemannian manifold is a smooth
$n-$dimensional manifold $M$ with a $k-$dimensional subbundle or
distribution $\triangle \subset TM$ on which a smooth inner product
$g_c$ is endowed. $(\triangle, g_c)$ is called a subriemannian
structure on $M$ and $\triangle$  horizontal bundle. In this paper,
we always assume $M$ is connected and $\triangle$ satisfies the
so-called Chow-H\"{o}mander condition  which means that vector
fields of $\triangle$ together with all their commutators span the
tangent space at each point on $M$. Carnot groups are important
examples of subriemannian manifolds. A Carnot group $\mathbb{G}$ is
a connected, simply connected Lie group with a graded Lie algebra
\begin{equation}\label{structure}
\flat=V^1\oplus\cdots\oplus V^r,\textrm{ with }
V^i=[V^1,V^{i-1}],[V^1,V^r]=0, i=2,\cdots,r.
\end{equation}
The integer $r$ is called the step of $\mathbb{G}$. Since
$\mathbb{G}$ is connected and simply connected, the exponential map
from $\flat$ to $\mathbb{G}$ gives  a global chart for $\mathbb{G}$.
Carnot groups are tangent spaces (in the sense of Gromov-Hausdorff)
of equiregular subriemannian manifolds, see \cite{Mi, Be}. It is
believed that the role played by Carnot groups in subriemannian
geometry is similar to the role of Euclidean Spaces in Riemannian
geometry.

 It follows from the Chow-Rashevskii connectivity theorem that
for any given points $p,q\in M$ there always exists at least a
horizontal curve connecting $p$ and $q$, see \cite{Chow,Ra}. Here a
horizontal curve is by definition an absolutely continuous curve
$\gamma:[0,1]\rightarrow M$ such that $\dot{\gamma}(t)\in
\Delta_{\gamma(t)}M$ whenever $\dot{\gamma}(t)$ exists. Thus one can
define a natural distance:
$${\rm d}_{sr}(p,q)=\inf\int_0^1\sqrt{g_c(\dot{\gamma},\dot{\gamma})}dt$$
where the infimum is taken among the set $\Omega(p,q)$ of all
horizontal curves $\gamma$ such that $\gamma(0)=p$ and
$\gamma(1)=q$. ${\rm d}_{sr}$ is called the
 Carnot-Carath\'{e}odory distance of $(M,\triangle,g_c)$. A subriemannian geodesic
 is a horizontal curve locally realizing ${\rm d}_{sr}$. We will reserve the terminology ``minimizing
 geodesic"
  or ``minimizer" for those globally distance-realizing subriemannian geodesics. It is not difficult to
  prove that any two sufficiently close points can be joined by a minimizing geodesic.
  If $(M,{\rm d}_{sr})$ is complete, there is a minimizing geodesic connecting any two points.
  Before Montgomery \cite{Mo2} (in 1991) discovered a smooth singular minimizer,
  it was taken for grant (see e.g. \cite{Str}) that each subriemannaian geodesic similar
  to a Riemannian geodesic could satisfy a Hamilton-Jacobi equation:
\begin{equation}\label{normal}
\dot{x}^i=\frac{\partial H}{\partial \lambda_i}, \qquad
\dot{\lambda}_i=-\frac{\partial H}{\partial x^i}
\end{equation}
where $(x^i,\lambda_i)$ is a coordinate system of $T^{\star}M$,
$H(x,\lambda)=\max_{v\in{\triangle_x}}\{\lambda(v)-\frac{1}{2}g_c(v,v)\}$
$(\lambda\in{T_x^\star M})$ is the subriemannian Hamiltonian. A
horizontal curve $\gamma(t)=(x^i(t))$ (denoted by a local
coordinate) satisfying (\ref{normal}) almost everywhere  for some
lift $\lambda(t)=(\lambda_i(t))$ can be proved to be locally
minimizing and smooth, and is called a normal geodesic. Montgomery's
example shows that not all subriemannian geodesics are normal. The
subriemannian geodesic problem is a special case of geometric
control problems. In fact singular curves or abnormal extremals play
a very important role  in optimal control theory, see e.g.
\cite{BC,AS0}.  It is well known that the Pontryagin Maximum
Principle (or the Lagrange Multiplier Rule in the Lagrangian
formulation) gives the first order necessary condition of optimality
for optimal control problems.  This first order condition is hardly
considered to be satisfactory when one studies abnormal extremals.
Recently experts developed necessary/sufficient second order
conditions of optimality, i.e, Goh condition and generalized
Legendre-Jacobi condition, see e.g. \cite{AG,AS0,AS1,AS2,AS3}. These
conditions were derived from the finiteness of the generalized Morse
index of critical points of the end-point mapping. We refer to
\cite{AS1,AS3,AS0} for the finiteness of the generalized Morse
index.

Let $\Omega(p)$ be the set of all horizontal curves
$\gamma:[0,1]\rightarrow M$ such that $\gamma(0)=p$. It is well
known that $\Omega(p)$ is a smooth Banach manifold, see
\cite{Bismut}. The end-point mapping is the smooth mapping
$\mathscr{E}:\Omega(p)\rightarrow M$ defined by taking
$\gamma\in{\Omega(p)}$ to $\gamma(1)$. Thus
$\Omega(p,q)=\mathscr{E}^{-1}(q)$. In general $\mathscr{E}$ is not
regular at all curves in $\Omega(p)$ and thus $\Omega(p,q)$ possibly
has no smooth structure if $q$ is a critical value of $\mathscr{E}$.
If $\gamma\in{\Omega(p)}$ is a critical point of $\mathscr{E}$, we
call $\gamma$ a singular curve. After Montgomery's example,
Liu-Sussmann in \cite{LS}  gave  more examples of singular curves
which are locally minimizing. All these examples found on rank two
distributions are in fact $C^1$-rigid curves which by definition are
locally isolated curves in $\Omega(p,q)$ with respect to the
$C^1$-topology, see also Bryant-Hsu \cite{BH}. For every rank two
distribution satisfying $\triangle_p^3\neq\triangle^2_p$ at $p\in M$
they proved that there exists a rigid curve locally
length-minimizing and emanating from $p$, see also Agrachev-Sarychev
in \cite[Theorem 6.2]{AS1}. Here $\triangle^1=\triangle,
\triangle^i=\triangle^{i-1}+[\triangle^1,\triangle^{i-1}]$ for
$i=2,\cdots$.  The research of such curves first appeared in the
work of Carath\'{e}odory, Engel, and Hilbert, see \cite{Bli,Young}.
Classical calculus of variations can not fully deal with the
subriemannian geodesic problem when $\Omega(p,q)$ contains singular
curves, because there possibly exist no smooth variations of such
curves. But singular curves could be subriemannian geodesics as
shown by the above mentioned work. A fundamental problem is whether
all subriemannian geodesics are smooth. This problem is equivalent
to the question whether all singular geodesics are smooth, since
normal (regular) geodesics are always smooth.  There are some
substantial results so far, while the problem is still open for
general cases. Gol\'{e}-Karidi in \cite{GK} constructed the first
example of a Carnot group  with a strictly singular minimizer and
they proved that all subriemannian geodesics in a step-two Carnot
group are normal. Agrachev-Sarychev \cite{AS2} proved that there
admit no strictly singular geodesics for medium fat distributions
including strong-generating distributions (fat distributions) for
which Strichartz \cite{Str} had already obtained the conclusion. For
a class of equiregular subriemannian manifolds, Leonardi-Monti
\cite{GP} showed that length-minimizing curves have no corner-like
singularities which in particular implies that all singular
geodesics in Carnot groups of rank two with step$\leq4$ are smooth.
There are also some ``generic" results which claims that for $3\leq
k<n$ there exists an open dense subset $\mathcal{O}_k$ of the space
$\mathscr{D}_k$ consisting of all k-dimensional distributions on $M$
(endowed with the Whitney $C^\infty$ topology), such that each
distribution in $\mathcal {O}_k$ admits no singular geodesics, see
\cite{AG, CJT} and references therein.

 In this paper we will concentrate on the case of Carnot groups.
 As mentioned above Carnot groups are very important in subriemannian geometry.
 Our study here will be instructive
 for later considerations of subriemannian geodesics of general distributions.
 Our starting point is the refined formulation
of the end-point mapping which boils down to the graded structures
 of Carnot groups. The graded structure (\ref{structure}) implies that
 each horizontal curve
 $\Upsilon=(\gamma^1,\gamma^2,\cdots,\gamma^r)$ is uniquely
 determined by the first layer $\gamma^1$. Here we use the
 exponential mapping $\exp$ to identify the Carnot group $\mathbb{G}$ with
 its Lie algebra $\flat$ and $\gamma^i=\pi^i(\exp^{-1}\Upsilon)\in
 V^i$ where $\pi^i:\flat=V^1\oplus\cdots\oplus V^r\rightarrow V^i$ is the
 projection to the $i$-th component.
The subriemannian geodesic problem in Carnot groups can be
formulated as a minimization problem with equality constraint. The
end-point mapping $\mathscr{E}$ is different from the ordinary one
which usually takes a control function to the end point of the
trajectory uniquely determined by the control function (the initial
point is fixed). The formula for the differential of $\mathscr{E}$
can be written out for Carnot groups with any step by a very tedious
computation. In the case of step 3, the differential and the
intrinsic Hessian of the end-point mapping are of simple form. We
will get the second order necessary conditions for optimality of a
singular curve, that is, if a singular curve
$\Upsilon=(\gamma^1,\gamma^2,\cdots,\gamma^r)$ is locally
energy-minimizing then $\gamma^1$ must satisfy the Goh condition and
the generalized Legendre-Jacobi condition, see Proposition \ref{goh}
and \ref{GLJ}. From the Goh condition and the graded structure
(\ref{structure}) we deduce that the first layer $\gamma^1$ of a
singular geodesic $\Upsilon$ must be in a lower-dimension subspace.
Thus singular geodesics must be in Carnot subgroups of rank 2 or
step 2. We finally reduce the problem to the rank two case which
(known to experts) is easy, see Theorem \ref{rank2} and Theorem
\ref{general}. The reduction to lower subgroups is a curious
coincidence with Hamenst\"{a}dt's suggestion for the smoothness
problem, see \cite{Ha}.

The paper is organized into five sections. In the next section we
give a formulation of the end-point mapping which is based on a
characterization of horizontal curves in Carnot groups. Section 3 is
devoted to the calculus of the end-point mapping. We will give the
differential, the intrinsic Hessian for the step$\leq3$ case. In
Section 4 we derive the second order necessary conditions of
singular geodesics. We prove the main results in Section 5.

\section{Horizontal Curves and the End-Point Mapping}
\subsection{Basic Structure of Carnot groups}
Let $\mathbb{G}$ be a Carnot group with a Lie algebra $\flat$
satisfying (\ref{structure}) (we call such Lie algebras Carnot
algebras). We endow on $V^1$ an inner product $<\cdot,\cdot>$. Let
$n_i=\dim(V_i)$, $n=\sum_{i=1}^rn_i$. The $r-$vector
$(n_1,n_1+n_2,\cdots,\Sigma_{j=1}^in_j,\cdots,n)$ is called the
growth vector of the Carnot group. We fix an orthonormal basis of
$V^1$, $\{e_1,\cdots,e_{n_1}\}$, then arbitrarily extend it to a
basis of $\flat$, $\{e_1,\cdots,e_{n_1},e_{n_1+1}\cdots,e_n\}$, and
extend $<\cdot,\cdot>$ to an inner product $g$ on $\flat$ making the
basis orthonormal.  Via the exponential mapping $\exp$ we identify
$\mathbb{G}$ with $\flat$ or $\mathbb{R}^n$ with a group law
determined by the Baker-Campbell-Hausdorff formula. For
$p\in{\mathbb{G}}$, setting
$X_i(p)=\left.\frac{d}{dt}\right|_{t=0}\{p\cdot\exp(te_i)\}$,
$i=1\cdots,n$, we get the basis of the space of left-invariant
vector fields. The horizontal bundle
$\triangle=span\{X_1,\cdots,X_{n_1}\}$ satisfies the
Chow-H\"{o}rmander condition by (\ref{structure}). Let
$g_c((X_i(p),X_j(p))=<e_i,e_j>, i,j=1,\cdots,n_1$. Thus we have a
subriemannian structure $(\triangle,g_c)$ on $\mathbb{G}$. We also
extend $g_c$ to a left-invariant Riemannian metric $g_r$ such that
$\{X_1,\cdots,X_n\}$ is an orthonormal basis of $T\mathbb{G}$. We
emphasize that subriemannian geodesics in $\mathbb{G}$ are
independent of the choice of orthonormal bases and their extensions,
that is, they are completely determined by $(\flat,
V^1,<\cdot,\cdot>)$ or equivalently by $(\mathbb{G},\triangle,g_c)$.

\begin{example} 1. The simplest Carnot group is the Heisenberg group $\mathbb{H}^m$ with
the Heisenberg algebra (growth vector=$(2m,2m+1)$)
$\mathbbm{h}=\textrm{span}\{e_1,\cdots,e_m,f_1\cdots,f_m\}\oplus
\textrm{span} \{g\}$ with the basis satisfying that $[e_i,f_i]=g$,
$i=1,\cdots,m$, and all other brackets vanish.

2. The Engel group is a Carnot group with the growth vector
$(2,3,4)$. Its algebra is
$\textrm{span}\{e_1,e_2\}\oplus\textrm{span}\{e_3\}\oplus\textrm{span}\{e_4\}$
with $[e_1,e_2]=e_3, [e_1,e_3]=e_4$. Note that  Carnot groups of
rank two $(n_1=2)$ has a special feature: the second layer has
dimension 1 whatever its step. We will see later this feature make
the study of its subriemannian geodesics very easy.

3. The free Carnot group with bi-dimension $(k,r)$ has the maximal
vector growth among all Carnot groups with $k$ generators and
step=$r$. The free Carnot groups play some particular roles in
nilpotent analysis.

4. Let $\flat$ be a Carnot algebra satisfying (\ref{structure}). If
$W\subset V^1$ is a lower-dimensional subspace, then
$\bar{\flat}=W^1\oplus\cdots\oplus W^{\bar{r}}$ is a Lie subalgebra
of $\flat$, where $W^1=W, W^{i}=[W^1,W^{i-1}]$,
$i=2,\cdots,\bar{r}$, and $\bar{r}\in[1,r]$ is the largest integer
such that $W^{\bar{r}}\neq0$. It is obvious that $\bar{\flat}$ is a
Carnot algebra and $\bar{\mathbb{G}}:=\exp(\bar{\flat})$ is a Carnot
subgroup of $\mathbb{G}$ (we regard Euclidean spaces as abelian
Carnot groups). We use $\bar{\flat}(W)$ (resp.
$\bar{\mathbb{G}}(W)$) to indicate the Carnot subalgebra (resp.
Carnot subgroup) generated by $W$. The reduction to Carnot subgroups
with lower-dimensional first layer is one of main tricks in this
paper.

\end{example}

Recall that the differential of the exponential mapping is given by
the following formula
\begin{equation}\label{diff}
    {\rm d}\exp(e)=\textrm{Id}-\sum_{m=2}^r\frac{(-1)^m}{m!}\textrm{ad}(e)^{m-1},
\end{equation}
see e.g. \cite{Va}. For an absolutely continuous curve
$\Upsilon:[0,1]\rightarrow \mathbb{G}$, we denote by $\gamma$ the
corresponding curve $\exp^{-1}(\Upsilon)$ with values in the Lie
algebra $\flat$. We have $\gamma=\Sigma_{i=1}^r\gamma^i$ where
$\gamma^i=\pi^i(\gamma)$, $\pi^i:\flat\rightarrow V^i$ is the
projection to the $i$-th layer. It is obvious that $\gamma$ is also
absolutely continuous. From (\ref{diff}) we get for a.e.
$t\in{[0,1]}$

\begin{equation} \label{derivative}
\begin{split}
\dot{\Upsilon}(t)&=\dot{\gamma}(t)-\sum_{m=2}^r\frac{(-1)^m}{m!}\textrm{ad}(\gamma(t))^{m-1}(\dot{\gamma}(t)) \\
                 &=\dot{\gamma}(t)-\sum_{m=2}^r\frac{(-1)^m}{m!}[\gamma(t),\dot{\gamma}(t)]_{m-1}
 \end{split}
 \end{equation}
In the last formula we used the iterated Lie bracket which is
defined by
$$
[e,f]_m=\underbrace{[e,[e,[\cdots,[e}_{m\textrm{
times}},f],],\cdots,] \textrm{ and }[e,f]_0=f.
$$
From (\ref{derivative}) we obtain that $\Upsilon$ is horizontal if
and only if for a.e. $t\in[0,1]$
$$
\pi^i\left(\dot{\gamma}(t)-\sum_{m=2}^r\frac{(-1)^m}{m!}[\gamma(t),\dot{\gamma}(t)]_{m-1}\right)=0,\quad
i=2,\cdots,r.
$$
We summarize as
\begin{lemma}\label{horizontal}
An absolutely curve $\Upsilon$ in $\mathbb{G}$ is horizontal if and
only if for a.e. $t\in[0,1]$
$$
\dot{\gamma}^i(t)=\sum_{m=2}^r\frac{(-1)^m}{m!}\pi^i([\gamma(t),\dot{\gamma}(t)]_{m-1}),
i=2,\cdots,r.
$$
\end{lemma}
We denote by $\mathcal{H}^1$ the Sobolev type space of all
horizontal curves $\Upsilon:[0,1]\rightarrow\mathbb{G}$  with square
integrable derivatives. In the rest of the paper we assume all
horizontal curves in $\mathbb{G}$ are  in $\mathcal{H}^1$. For our
purpose this assumption is not restrictive since all rectifiable
curves can be arc-length parameterized. Combining (\ref{derivative})
with Lemma \ref{horizontal}, we have for a.e. $t\in[0,1]$
\begin{equation}\label{submersion1}
\dot{\Upsilon}(t)=\sum_{i=1}^{n_1}\dot{x}^i(t)X_i(\Upsilon(t))
\end{equation}
where $\gamma^1(t)=\Sigma_{i=1}^{n_1}x^i(t)e_i$.  Define $\mathcal
{P}:\mathbb{G}\rightarrow V^1$, $\mathcal
{P}(p)=\pi^1(\exp^{-1}(p))$. The mapping $\mathcal {P}$ is just the
projection
$\mathbb{R}^n\ni(x_1,\cdots,x_{n_1},\cdots,x_n)\rightarrow
(x_1,\cdots,x_{n_1})\in\mathbb{R}^{n_1}$ when we identify
$\mathbb{G}$ as $(\mathbb{R}^n,\cdot)$. The formula
(\ref{submersion1}) in particular implies that $\mathcal {P}$ is a
Riemannian submersion from $(\mathbb{G},g_r)$ to
$(V^1,<\cdot,\cdot>)$ (or equivalently from $(\mathbb{R}^n,g_r)$ to
$(\mathbb{R}^{n_1},<\cdot,\cdot>)$) with  the property that for any
$p\in G$, $\mathcal {P}_{\star,p}(X_i(p))=e_i$, $i=1,\cdots,n_1$,
and $\mathcal {P}_{\star,p}(X_j(p))=0$ for $j=n_1+1,\cdots,n$.

Note that the graded condition (\ref{structure}) for the Lie algebra
$\flat$ is equivalent to the following condition
\begin{equation}\label{structure2}
    V^i=\underbrace{[V^1,[\cdots,[V^1,V^1]]]}_{i\textrm{ times}},\quad
    i=2,...r,\textrm{ and }V^j=0\textrm{ for }j>r,
\end{equation}
which together with Lemma \ref{horizontal} implies that
\begin{equation}\label{recursive}
\dot{\gamma}^i=\sum_{m=2}^{i}\frac{(-1)^m}{m!}\left(\sum_{j_1+j_2+\cdots+j_m=i}
[\gamma^{j_1},[\gamma^{j_2},[\cdots,[\gamma^{j_{m-1}},\dot{\gamma}^{j_m}]]]]\right)
\textrm{ for }i=2,\cdots,r,\textrm{ a.e..}
\end{equation}
which means that $\gamma^1$ determines
$\gamma^2,\gamma^3,\cdots,\gamma^r$ recursively. We list
$\dot{\gamma}^2,\dot{\gamma}^3,\dot{\gamma}^4$ as functions of
$\gamma^1$:

\begin{equation} \label{determine}
\left\{ \begin{aligned}
         \dot{\gamma}^2&= \frac{1}{2}[\gamma^1,\dot{\gamma}^1] \\
         \dot{\gamma}^3&=\frac{1}{2}\left\{[\gamma^1,\dot{\gamma}^2]+
         [\gamma^2,\dot{\gamma}^1]\right\}-\frac{1}{6}[\gamma^1,[\gamma^1,\dot{\gamma}^1]]\\
         \dot{\gamma}^4&=\frac{1}{2}\left\{[\gamma^1,\dot{\gamma}^3]+[\gamma^2,\dot{\gamma}^2]+[\gamma^3,\dot{\gamma}^1]\right\}
         -\frac{1}{6}\left\{[\gamma^1,[\gamma^1,\dot{\gamma}^2]]+[\gamma^1,[\gamma^2,\dot{\gamma}^1]]\right.\\&\quad+
         \left.[\gamma^2,[\gamma^1,\dot{\gamma}^1]]\right\}+\frac{1}{24}[\gamma^1,\dot{\gamma}^1]_3.
                          \end{aligned} \right.
                          \end{equation}

Sometimes we will abuse the notation
$\Upsilon=(\gamma^1,\cdots,\gamma^r)$ or
$\Upsilon=\Sigma_{i=1}^r\gamma^i$.
\begin{prop}\label{lift}{\rm(1)} Given $p\in{\mathbb{G}}$. Every absolutely continuous curve
 $\gamma^1:[0,1]\rightarrow V^1$ has a unique horizontal lift
$\Upsilon=(\gamma^1,\cdots,\gamma^r):[0,1]\rightarrow\mathbb{G}$
determined by (\ref{recursive}) with $\Upsilon(0)=p$. They have the
same length and same regularity or smoothness;
     {\rm(2)} Horizontal lifts of every straight line (or its interval) in
$V^1$ are subriemannian minimizing geodesics.
\end{prop}
\proof (1) Note that the class of absolutely continuous curves in
$\flat$ is just the Sobolev class $W^{1,1}([0,1],\flat)$. This
implies $\gamma^1$ is continuous and thus bounded in $[0,1]$ with
$\dot{\gamma}^1\in L^1([0,1],\flat)$, see e.g. \cite[Chapter
2]{BGH}. So there exists $\gamma^2\in W^{1,1}([0,1],\flat)$ such
that
$\gamma^2(t)=\frac{1}{2}\int_0^t[\gamma^1,\dot{\gamma}^1]d\tau+\pi^2(\exp^{-1}p)$.
Continuing this process, and noting that each summand in the right
hand side of (\ref{recursive}) contains only one term with
derivative and other terms are bounded, that is, the right hand side
of (\ref{recursive}) is in $L^1([0,1],\flat)$. Thus the function
$\gamma^i$ satisfying
$$
\gamma^i(t)=\int_0^t\sum_{m=2}^{i}\frac{(-1)^m}{m!}\left(\sum_{j_1+j_2+\cdots+j_m=i}
[\gamma^{j_1},[\gamma^{j_2},[\cdots,[\gamma^{j_{m-1}},\dot{\gamma}^{j_m}]]]]\right)d\tau+\pi^i(\exp^{-1}p)
$$
is in $W^{1,1}([0,1],\flat)$. From (\ref{submersion1}) $\gamma^1$
and its lift above have the same length.

 To see (2), we recall by
definition that in $(\mathbb{G},g_r)$ Riemannian geodesics which are
horizontal must be subriemannian
 geodesics. Since the geodesics of Euclidean space
$(V^1,<\cdot,\cdot>)$ are straight lines (or their intervals) ,
their horizontal lifts are Riemannian geodesics because $\mathcal
{P}$ is a Riemannian submersion, see \cite{Oneil}. The minimizing
property is obvious.
\endproof

By Proposition \ref{lift} we sometimes do not distinguish a
horizontal curve $\Upsilon=(\gamma^1,\cdots,\gamma^r)$ with its
projection to the first layer $\gamma^1$.

Note that horizontal lifts of straight lines in $V^1$ are not
necessarily still a line (looking in $\flat$) if $\mathbb{G}$ with
step$\geq3$. In fact, let $\gamma^1(t)=vt+v_0$ with $v,v_0\in{V^1}$.
By the formula (\ref{determine}) we have
$$
\left\{ \begin{aligned}
         \dot{\gamma}^2(t)&= \frac{1}{2}[v,v_0]\\
                  \dot{\gamma}^3(t)&=\frac{1}{6}[v,[v,v_0]]t-\frac{5}{12}[v_0,[v_0,v]]+\frac{1}{4}[\gamma^2(0),v].
                          \end{aligned} \right.
$$
So the third layer of the lift is not a line unless $[v,[v,v_0]]=0$.
While if $v_0=0$, that is, the line passes through the origin, its
lift is just itself.

\subsection{The End-Point Mapping} Given $p,q\in{\mathbb{G}}$, deonte by
$\Omega(p)$ the Hilbert manifold of all horizontal curves
$\Upsilon\in\mathcal {H}^1$ with $\Upsilon(0)=p$. Let
$\Omega(p,q)=\{\Upsilon\in\Omega(p):\Upsilon(1)=q\}$.
 Since  $\Omega(p)=p\cdot\Omega(0):=\{p\cdot\Upsilon: \Upsilon\in\Omega(0)\}$, $\Omega(p,q)=p\cdot\Omega(0,p^{-1}\cdot q)$
 and the metric $g_c$ is left-invariant, it suffices to
consider horizontal curves emanating from the unit. Here we abuse 0
to denote the unit of $\mathbb{G}$. From Proposition \ref{lift}, we
see that the projection $\mathcal{P}$ gives a bijective mapping
(still denoted by $\mathcal {P}$) from $\Omega(0)$ to
$H^1(0):=\{\gamma^1\in H^1([0,1],V^1):\gamma^1(0)=0\}$ with the
mapping of horizontal lift as its inverse, where $H^1([0,1],V^1)$
denote the Sobolev space of all absolutely curves
$\gamma^1:[0,1]\rightarrow V^1$ with square integrable derivatives.
Note from the formula (\ref{recursive}), we have for
$\Upsilon\in{\Omega(0)},t\in[0,1]$,
$$
\Upsilon(t)=(\gamma^1(t),\gamma^2(t),\cdots,\gamma^r(t))
$$
where $\gamma^i(t), i=2,\cdots,r$ is regard as a mapping $F^{i,t}$
(defined recursively) from $H^1(0)$ to $V^i$:
\begin{equation}\label{dayu2}
F^{i,t}(\gamma^1)=\int_0^t\sum_{m=2}^{i}\frac{(-1)^m}{m!}\left(\sum_{j_1+j_2+\cdots+j_m=i}
[\gamma^{j_1},[\gamma^{j_2},[\cdots,[\gamma^{j_{m-1}},\dot{\gamma}^{j_m}]]]]\right)d\tau.
\end{equation}

 Now the original end-point mapping
\begin{equation} \label{e-p-m}
\textrm{end}:\Omega(0)\ni\Upsilon\rightarrow\Upsilon(1)\in\mathbb{G}
 \end{equation}
can be interpreted as
\begin{equation}\label{epm}
\mathscr{E}:H^1(0)\ni \gamma\rightarrow
(F^1(\gamma^1),F^2(\gamma^1),\cdots,F^r(\gamma^1))\in\flat
\end{equation}
where $F^{1}(\gamma^1)=\gamma^1(1),F^{i}=F^{i,1}, i=2,\cdots,r$.

Noting that given $\exp\xi=q\in\mathbb{G}$ for $\xi\in{\flat}$,
$\Omega(0,q)=\exp(\mathscr{E}^{-1}(\xi)),$  the subriemannian
geodesic problem in $\mathbb{G}$
$$
\min_{\Upsilon\in{\Omega(0,q)}}\frac{1}{2}\int_0^1g_c(\dot{\Upsilon},\dot{\Upsilon})dt
$$
is equivalent to the minimizing problem with equality constraint
\begin{equation}\label{problem}
\min_{\mathscr{E}(\gamma^1)=\xi}\frac{1}{2}\int_0^1\left\vert\dot{\gamma}^1\right\vert^2dt.
\end{equation}
By the Cauchy-Schwarz inequality  the problem of minimizing the
energy functional is equivalent to that of minimizing  the length
functional. The existence of the subriemannian geodesic problem even
for general subriemannian manifolds can be obtained by an argument
of direct method in calculus of variations, see e.g. Appendix D in
\cite{Mo1} where one also will find the ordinary formulation of the
end-point mapping. What we are concerned with is their smoothness.
The refined mapping $\mathscr{E}:H^1(0)\rightarrow \flat$ from a
Hilbert space to a vector space will help us much.

\section{The Calculus Of The End-Point Mapping}
\subsection{The Generalized Morse Index Theorem} To begin some computation, let us first see what we need
according to the theory of generalized Morse index which we will
later resort to. Let $(X,\parallel\cdot\parallel)$ be a Banach
space, $L:X\rightarrow \mathbb{R}$, and $F:X\rightarrow Y$ with $Y$
a finite dimensional vector space, be $C^2$ Fr\'{e}chet
differentiable mappings. Given $y_0\in{Y}$, consider the minimizing
problem with equality constraint
\begin{equation}\label{Lagrange1}
    \min_{F(x)=y_0} L(x).
\end{equation}
The  Lagrange Multiplier Rule states that if $x\in X$ is a solution
of (\ref{Lagrange1}) then there exists a nontrivial couple
$(\lambda^0,\lambda^\star)\in{\mathbb{R}\times Y^\star}$ such that
$\lambda^0{\rm d}_xL+\lambda^\star {\rm d}_xF=0$, where ${\rm d}_xL$
(resp. ${\rm d}_xF$) denotes the Fr\'{e}chet derivative of $L$
(resp. $F$) at the point $x$. In other words, the point $x$ is a
singular point of the augmented end-point mapping $\mathcal {L}:X\ni
y\rightarrow (L(y),F(y))\in \mathbb{R}\times Y $ and
\begin{equation}\label{lagrange3}
\tilde{\lambda}\cdot{\rm d}_x\mathcal{L}=\lambda^0{\rm
d}_xL+\lambda^\star {\rm d}_xF=0\textrm{ with
}\tilde{\lambda}=(\lambda^0,\lambda^\star).
\end{equation}
 The abnormal case $\lambda^0=0$ arises
exactly when Im(d$_xF)\neq Y$, i.e., $x$ is a singular point of $F$.
In this case we call $(x,\lambda^\star)$ an abnormal extremal. In
the regular case Im(d$_xF)= Y$,  we take $\lambda^0=1$. For $x\in X$
if there exists $\lambda^\star$ such that
\begin{equation}\label{lagrange2}
{\rm d}_xL+\lambda^\star {\rm d}_xF=0,
\end{equation}
we call $(x,\lambda^\star)$ a normal extremal.   The definition of
normal extremals is equivalent to the one given in the Introduction.
In the theory of subriemannian geodesics there exists a
correspondence between $\lambda^\star$ in (\ref{lagrange2}) and the
Hamiltonian lift $\lambda(t)$ in (\ref{normal}), see \cite{Hsu} or
\cite[Chapter 5]{Mo1}. We remark that an abnormal extremal may be
normal by choosing a suitable multiplier (or a Hamiltonian lift).
Those abnormal extremals which can not be normal for any multiplier
are called strictly abnormal extremals.

The corank of $x$ is defined as the codimension of Im(d$_x\mathcal
{L})$. The following theorem, which gives necessary/sufficient
conditions of optimality for a singular point is enough for our
purpose, for the general versions see \cite[Chapter 20]{AS0}.
\begin{theorem}[\cite{AG2,AS2,AS1,AS0}]\label{Morse} If $x$ is a local minimizer in $X$ of the minimizing
problem (\ref{Lagrange1}), of corank $N$, then for the nontrivial
pair of Lagrange multiplier
$\tilde{\lambda}=(\lambda^0,\lambda^\star)$ ($\lambda^0=0\textrm{ or
}1$) satisfying  (\ref{lagrange3}), the Morse index of the quadratic
form $\tilde{\lambda}\cdot{\rm d}^2_x\mathcal{L}$ restricted to
$\ker {\rm d}_xF$ is less than or equal to $N-1$.
\end{theorem}
 We recall that the Morse index of a quadratic form is
the maximal dimension of subspaces  on which the quadratic form is
negative definite. Theorem \ref{Morse} is classical for the regular
case for which (\ref{lagrange2}) is the Euler-Lagrange equation, see
\cite{Mil}.

\subsection{The Differential of The End-Point Mapping}
In the next section we will derive second order necessary conditions
for minimality of abnormal extremals of the problem (\ref{problem})
from the finiteness of the Morse index of the quadratic form
$\tilde{\lambda}\cdot{\rm d}^2_x\mathcal{L}$ stated in Theorem
\ref{Morse}. In the following we do some computation of the
differential of the end-point mapping.
\begin{lemma}
Given $\Upsilon=(\gamma^1,\cdots,\gamma^r)\in\Omega(0)$, then
$T_\Upsilon\Omega(0)=T_{\gamma^1}H^1(0)=H^1(0).$
\end{lemma}
\proof $H^1(0)$ is a Hilbert space. Let $\phi\in{H^1(0)}$. The
family of horizontal lifts of $\gamma^1+\epsilon\phi$
($\epsilon\in[-\epsilon_0,\epsilon_0])$,
$$
\Upsilon_\epsilon(t)=\left(\gamma^1(t)+\epsilon\phi(t),F^{2,t}(\gamma^1+\epsilon\phi),\cdots,F^{r,t}(\gamma^1+\epsilon\phi)\right)
$$
where $F^{i,t}, i=2\cdots,r,$ is defined as in (\ref{dayu2}), is a
smooth curve in $\Omega(0)$ with $\Upsilon_0=\Upsilon$.

On the other hand, if $\Upsilon_\epsilon$ is a smooth family in
$\Omega(0)$ with $\Upsilon_0=\Upsilon$, then by Lemma
\ref{horizontal} and (\ref{recursive}), (\ref{dayu2}) we have
$$
\Upsilon_\epsilon(t)=\left(\gamma_\epsilon^1,F^{2,t}(\gamma^1_\epsilon),\cdots,F^{r,t}(\gamma^1_\epsilon)\right)
$$
where $\gamma^1_\epsilon=\mathcal {P}(\Upsilon_\epsilon)$ is a
smooth family in $H^1(0)$ with
$\gamma^1_0=\gamma^1=\mathcal{P}(\Upsilon)$.
\endproof
For  $\gamma^1,\phi\in{H^1(0)}$, the differential at $\gamma^1$ of
the end-point mapping $\mathscr{E}$ is
\begin{equation}\label{differential}
{\rm d}_{\gamma^1}\mathscr{E}:H^1(0)\ni\phi\rightarrow{\rm
d}_{\gamma^1}\mathscr{E}(\phi)=\left(\phi(1),{\rm
d}_{\gamma^1}F^2(\phi),\cdots,{\rm
d}_{\gamma^1}F^r(\phi)\right)\in{T_q}\mathbb{G}
\end{equation}
where $q=(\gamma^1(1),F^2(\gamma^1),\cdots,F^r(\gamma^1))$ and $F^i$
is shortened for $F^{i,1}$, $i=2,\cdots,r$. Noting that the
differential ${\rm d}_{\gamma^1}F^i(\phi)$, $i=2,\cdots,r$,
recursively depends on
 ${\rm d}_{\gamma^1}\dot{F}^{j,t}(\phi)$, ${\rm d}_{\gamma^1}F^{j,t}(\phi),
j=2,\cdots,i-1,$ $t\in(0,1]$, its computation is complicated for
$i\geq5$.
 We restrict to the case
of step$\leq3$ partly also because of technical difficulties  in the
next two sections. Observe that $\frac{\rm d}{\rm dt}{\rm
d}_{\gamma^1}F^{i,t}(\phi)={\rm d}_{\gamma^1}\dot{F}^{i,t}(\phi)$.
From (\ref{determine})-(\ref{dayu2}) we have for step=2
\begin{equation} \label{F2}
\left\{ \begin{aligned}
         F^{2,t}(\gamma^1)&=\frac{1}{2}\int_0^t[\gamma^1,\dot{\gamma}^1]d\tau \\
         {\rm d}_{\gamma^1}F^{2,t}(\phi)
         &=\int_0^t[\phi,\dot{\gamma}^1]d\tau+\frac{1}{2}[\gamma^1(t),\phi(t)]\\
         {\rm
         d}_{\gamma^1}\dot{F}^{2,t}(\phi)&=\frac{1}{2}[\phi(t),\dot{\gamma}^1(t)]+\frac{1}{2}[\gamma(t),\dot{\phi}(t)]
                          \end{aligned} \right.
                          \end{equation}
For step=3, from (\ref{determine}) we first have
$$
\dot{F}^{3,t}=\dot{\gamma}^3(t)=\frac{1}{2}\frac{\rm d}{\rm
dt}[F^{2,t},\gamma^1(t)]+\frac{1}{3}[\gamma^1(t),[\gamma^1(t),\dot{\gamma}^1(t)]],
$$
then using (\ref{F2}) get
\begin{equation} \label{F3}
\left\{ \begin{aligned}
         {\rm d}_{\gamma^1}F^{3,t}(\phi)=&\int_0^t[\gamma^1,[\phi,\dot{\gamma}^1]]d\tau+
         \frac{1}{2}\left([F^{2,t},\phi]+[{\rm d}_{\gamma^1}F^{2,t}(\phi),\gamma^1]\right)+\frac{1}{3}[\gamma^1(t),\phi(t)]_2 \\
                  =&\int_0^t[\gamma^1,[\phi,\dot{\gamma}^1]]d\tau+\frac{1}{2}\left[\int_0^t[\phi,\dot{\gamma}^1]d\tau,\gamma^1(t)\right]
                  +\frac{1}{4}\left[\int_0^t[\gamma^1,\dot{\gamma}^1]d\tau,\phi(t)\right]\\
                  &+\frac{1}{12}[\gamma^1(t),[\gamma^1(t),\phi(t)]]
\end{aligned} \right.
\end{equation}
In the above computation we used integration by parts, skew-symmetry
 and Jacobi identity of Lie brackets to arrange terms.
\begin{lemma}
In the case $r=3$, by (\ref{differential})-(\ref{F3}) we have
$\phi\in\ker\left({\rm d}_{\gamma^1}\mathscr{E}\right)$ if and only
if
\begin{equation} \label{kernel}
    \left.\begin{aligned}
             \phi(1) &= 0 \\
                      \int_0^1[\phi,\dot{\gamma}^1]dt&=0\\
                      \int_0^1[\gamma^1,[\phi,\dot{\gamma}^1]]dt&=0
                              \end{aligned} \right\}
                              \end{equation}
\end{lemma}
Now we compute the second Fr\'{e}chet derivative of the end-point
mapping $\mathscr{E}$ for $r=3$. From (\ref{F2}) we have
\begin{equation}\label{2F2}
{\rm d}^2_{\gamma^1}F^2(\phi,\phi)=\int_0^1[\phi,\dot{\phi}]dt.
\end{equation}
For $\phi\in\ker\left({\rm d}_{\gamma^1}\mathscr{E}\right)$ it
follows from (\ref{F3}) and (\ref{kernel}) that
\begin{equation}\label{2F3}
{\rm
d}^2_{\gamma^1}F^3(\phi,\phi)=\int^1_0[\phi,[\phi,\dot{\gamma}^1]]dt+\int_0^1\left[\gamma^1-\frac{1}{2}\gamma^1(1),[\phi,\dot{\phi}]\right]dt.
\end{equation}
So the intrinsic quadratic mapping (see e.g. \cite[p.294-296]{AS0})
of $\mathscr{E}$ for $r=3$
$$
{\rm d}^2_{\gamma^1}\mathscr{E}:\ker\left({\rm
    d}_{\gamma^1}\mathscr{E}\right)\times\ker\left({\rm
    d}_{\gamma^1}\mathscr{E}\right)\rightarrow\flat
$$
is
\begin{equation}\label{2epm}
{\rm
d}^2_{\gamma^1}\mathscr{E}(\phi,\phi)=\left(0,\int_0^1[\phi,\dot{\phi}]dt,
\int^1_0[\phi,[\phi,\dot{\gamma}^1]]dt+\int_0^1
\left[\gamma^1-\frac{1}{2}\gamma^1(1),[\phi,\dot{\phi}]\right]dt\right).
\end{equation}

\section{Goh Condition and Legendre-Jacobi Condition}
In the rest of the paper we assume the Carnot group has step$\leq3$.
From (\ref{F2}) and (\ref{F3}) we have for $\phi\in H^1(0)$ with
$\phi(1)=0$
\begin{equation}\label{simple1}
{\rm
d}_{\gamma^1}\mathscr{E}(\phi)=\left(0,\int_0^1[\phi,\dot{\gamma}^1]dt,
            \int_0^1\left[\gamma^1
-\frac{1}{2}\gamma^1(1),\left[\phi,\dot{\gamma}^1\right]\right]dt\right).
\end{equation}

Applying the Lagrange Multiplier Rule to the problem (\ref{problem})
we get
\begin{prop} If $\gamma^1\in{H^1(0)}$ is a minimizer of the problem
(\ref{problem}), then there exists
$\tilde{\lambda}=(\lambda^1,\lambda^2,\lambda^3)\in\flat^\star$ with
$\lambda^i\in{(V^i)^\star}$, $i=1,2,3$, such that for any $\phi\in
H^1(0)$
\begin{equation}\label{normal1}
\int_0^1<\dot{\gamma}^1,\dot{\phi}>+\lambda^2\left[\phi,\dot{\gamma}^1\right]+\lambda^3\left[\gamma^1
-\frac{1}{2}\gamma^1(1),\left[\phi,\dot{\gamma}^1\right]\right]dt=0
\textrm{ with }\phi(1)=0
\end{equation}
or
\begin{equation}\label{singular}
\tilde{\lambda}{\rm
d}_{\gamma^1}\mathscr{E}(\phi)=\lambda^1\phi(1)+\lambda^2{\rm
d}_{\gamma^1}F^2(\phi)+\lambda^3 {\rm d}_{\gamma^1}F^3(\phi)=0
\textrm{ with }\tilde{\lambda}\neq0.
\end{equation}
\end{prop}
\proof Taking $L$ as the energy functional
$L(\gamma^1)=\frac{1}{2}\int_0^1\vert\dot{\gamma}^1\vert^2dt$, by
the Lagrange Multiplier Rule there exists nontrivial
$(\lambda^0,\tilde{\lambda})\in\mathbb{R}\times\flat^\star$ with
$\lambda^0=0\textrm{ or }1$ and
$\tilde{\lambda}=(\lambda^1,\lambda^2,\lambda^3)$ such that for any
$\phi\in H(0)$
$$
\lambda^0{\rm d}_{\gamma^1}L(\phi)+\lambda^1\phi(1)+\lambda^2{\rm
d}_{\gamma^1}F^2(\phi)+\lambda^3{\rm d}_{\gamma^1}F^3(\phi)=0.
$$
When $\lambda^0=1$, (\ref{normal1}) follows from the last formula,
(\ref{simple1}) and ${\rm
d}_{\gamma^1}L(\phi)=\int_0^1<\dot{\gamma}^1,\dot{\phi}>dt$.
\endproof

We call $\gamma^1$ (or its horizontal lift) satisfying
(\ref{normal1}) for some $(\lambda^2,\lambda^3)$ a normal geodesic.
By a standard argument from the theory of (elliptic) differential
equations normal geodesics are smooth, see e.g. \cite{BGH}. Singular
geodesics are local minimizers $\gamma^1$ (or their horizontal
lifts) of the problem (\ref{problem}) satisfying (\ref{singular})
for some
$\tilde{\lambda}=(\lambda^1,\lambda^2,\lambda^3)\in\flat^\star$. For
singular geodesics the following result follows from Theorem
\ref{Morse} and (\ref{2epm}).
\begin{prop}\label{quadratic1}
If $\gamma^1\in H^1(0)$ is a singular minimizer of the problem
(\ref{problem}), there exists a nontrivial
$\tilde{\lambda}=(\lambda^1,\lambda^2,\lambda^3)$ satisfying
(\ref{singular}) such that the Morse index of the quadratic form
\begin{equation}\label{quadratic}
\begin{aligned}
\tilde{\lambda}{\rm
d}^2_{\gamma^1}\mathscr{E}(\phi,\phi)=&\int_0^1\lambda^2[\phi,\dot{\phi}]dt+\int_0^1
\lambda^3\left[\gamma^1-\frac{1}{2}\gamma^1(1),[\phi,\dot{\phi}]\right]dt\\
&+\int^1_0\lambda^3[\phi,[\phi,\dot{\gamma}^1]]dt\qquad
\left(\phi\in{\ker\left({\rm d}_{\gamma^1}\mathscr{E}\right)}\right)
\end{aligned}
\end{equation}
is finite.
\end{prop}
\begin{prop}\label{goh}
Let $\gamma^1,\tilde{\lambda}=(\lambda^1,\lambda^2,\lambda^3)$ be as
in Proposition \ref{quadratic1}. Assume $\gamma^1$ is parameterized
proportionally to arc-length. Then
\begin{equation} \label{eq:goh}
\left\{ \begin{aligned}
         \lambda^2 &= 0 \\
                  \lambda^3\left[\gamma^1(t),[a,b]\right]&=0,\quad \forall
                  a,b\in{V^1},\forall t\in[0,1].
                          \end{aligned} \right.
                          \end{equation}
\end{prop}
\proof The argument is similar to \cite[Proposition 20.13]{AS0}. The
idea is a type of   scaling or blowing up method.

Let $\bar{\tau}\in[0,1]$ be a Lebesgue point of $\dot{\gamma}^1$.
Take a smooth mapping $c:\mathbb{R}\rightarrow V^1$ with support on
$[0,2\pi]$ such that $\int_0^{2\pi}c(s)ds=0$. Let
$\phi_{c,\epsilon}(\tau)=\int_0^\tau
c(\frac{\tilde{\tau}-\bar{\tau}}{\epsilon})d\tilde{\tau}$ with
$\epsilon$ small. This certainly implies
$\dot{\phi}_{c,\epsilon}(\tau)=c(\frac{\tau-\bar{\tau}}{\epsilon})$
for  $\tau\in[0,1]$ and $\phi_{c,\epsilon}\in{H^1(0)}$ with
$\phi_{c,\epsilon}(1)=0$. Letting $w(s)=\int_0^s
c(\tilde{s})d\tilde{s}$, we have
\begin{equation}\label{goh-eq1}
\int_0^1\lambda^2[\phi_{c,\epsilon},\dot{\phi}_{c,\epsilon}]d\tau=\epsilon^2\int_0^{2\pi}\lambda^2\left[w(s),c(s)\right]ds
\end{equation}
 and similarly
\begin{equation}\label{goh-eq2}
\begin{aligned}
\int_0^1\lambda^3\left[\gamma^1-\frac{1}{2}\gamma^1(1),[\phi_{c,\epsilon},\dot{\phi}_{c,\epsilon}]\right]dt=&\epsilon^2\int_0^{2\pi}\lambda^3\left[\gamma^1(\bar{\tau})
-\frac{1}{2}\gamma^1(1), \left[w(s),c(s)\right]\right]ds\\
&+\epsilon^3O(1)
\end{aligned}
\end{equation}
where we used the fact $\dot{\gamma}\in L^\infty$ and thus
$\gamma^1(\bar{\tau}+\epsilon s)=\gamma^1(\bar{\tau})+\epsilon O(1)$
for $\epsilon$ small enough. For the last term in (\ref{quadratic})
we have
\begin{equation}\label{goh-eq3}
\int^1_0\lambda^3[\phi_{c,\epsilon},[\phi_{c,\epsilon},\dot{\gamma}^1]]dt=\epsilon^3\int_0^{2\pi}\lambda^3\left[w(s),\left[w(s),\dot{\gamma}^1(\bar{\tau}+\epsilon
s)\right]\right]ds=\epsilon^3O(1).
\end{equation}
From (\ref{goh-eq1})-(\ref{goh-eq3}), we get
\begin{equation}\label{goh-eq4}
\tilde{\lambda}{\rm
d}^2_{\gamma^1}\mathscr{E}(\phi_{c,\epsilon},\phi_{c,\epsilon})=\epsilon^2\int_0^{2\pi}\omega(w(s),c(s))ds+\epsilon^3O(1)
\end{equation}
where
$\omega(a,b)=\lambda^2[a,b]+\lambda^3\left[\gamma^1(\bar{\tau})
-\frac{1}{2}\gamma^1(1), [a,b]\right]$ is a skew-symmetric bilinear
form on $V^1$.

We claim that $\omega\equiv0$. In fact, if $\omega\neq0$, then ${\rm
rank}\omega=2l_0>0$ and we can change the basis of $V^1$ such that
$$
\omega(a,b)=\sum_{i=1}^{l_0}\left(x^iy^{i+l_0}-x^{i+l_0}y^i\right)
$$
for any $b=(x^1,\cdots,x^{n_1}),a=(y^1,\cdots,y^{n_1})\in V^1$. Now
we take
$$
c(s)=\left(x^1(s),0,\cdots,0,x^{l_0+1}(s),0,\cdots,0\right)
$$
where $ x^1(s)=\sum^\infty_{k=1}\xi_k\cos ks,
x^{l_0+1}(s)=\sum_{k=1}^\infty\eta_k\sin ks$, and
$(\xi_k)_{k=1}^\infty,(\eta_k)_{k=1}^\infty\in l^1$. Putting
$c(s),w(s)=\int_0^sc(\tilde{s})d\tilde{s}$ into (\ref{goh-eq4}) we
get
$$
\tilde{\lambda}{\rm
d}^2_{\gamma^1}\mathscr{E}(\phi_{c,\epsilon},\phi_{c,\epsilon})=-\left(2\pi\sum_{k=1}^\infty\frac{1}{k}\xi_k\eta_k\right)\epsilon^2+\epsilon^3O(1).
$$
From the last formula and the construction of $\phi_{c,\epsilon}$ it
follows that there exists an infinite dimensional space $\mathcal
{K}$ such that $\tilde{\lambda}{\rm
d}^2_{\gamma^1}\mathscr{E}(\phi,\phi)<0$ for each $\phi\in\mathcal
{K}$. Note that $\mathcal {K}\cap \ker\left({\rm
d}_{\gamma^1}\mathscr{E}\right)$ is also infinite dimensional, since
the rank of $\mathscr{E}$ is less than $n$. It implies the Morse
index of $\tilde{\lambda}{\rm d}^2_{\gamma^1}\mathscr{E}$ is
infinite. This is impossible by Proposition \ref{quadratic1}, so
$\omega\equiv0$.

We have proved that if $t\in{[0,1]}$ is a Lebesgue point of
$\dot{\gamma}^1$, then
\begin{equation}\label{goh-eq6}
\lambda^2[a,b]+\lambda^3\left[\gamma^1(t) -\frac{1}{2}\gamma^1(1),
[a,b]\right]=0 \textrm{ for any }a,b\in{V^1}.
\end{equation}
Since almost all points in $[0,1]$ are Lebesgue points of
$\dot{\gamma}^1\in L^\infty$ by the Lebesgue differentiation
theorem, it follows from the continuity of $\gamma^1$ that
(\ref{goh-eq6}) holds for any $t\in[0,1]$. In (\ref{goh-eq6})
letting $t=0$ we get $\lambda^2[a,b]+\lambda^3\left[
-\frac{1}{2}\gamma^1(1), [a,b]\right]=0$ (since $\gamma^1(0)=0$).
Combing the last identity with (\ref{goh-eq6}), we obtain
$\lambda^3\left[\gamma^1(t), [a,b]\right]=0$ for any $t\in[0,1]$ and
any $a,b\in V^1$. Applying the identity (\ref{goh-eq6}) again, we
finally have $\lambda^2\xi=0$ for any $\xi\in V^2$, since
$[V^1,V^1]=V^2$.
\endproof

The condition of (\ref{eq:goh}) is called the Goh condition which
first appeared in references of singular control theory, see
\cite{Goh}. A curve $\gamma^1\in{H^1(0)}$ (or its horizontal lift)
together with some nonzero
$\tilde{\lambda}=(\lambda^1,\lambda^2,\lambda^3)\in\flat^\star$
satisfying the Goh condition is called a Goh curve and the pair
$(\gamma^1,\tilde{\lambda})$ is called a Goh extremal. Note that for
a Goh extremal $(\gamma^1,\tilde{\lambda})$, $\lambda^3\neq0$ and
(\ref{singular}) automatically holds by choosing $\lambda^1=0$. The
following fact is instructive for the study of subriemannian
geodesics even for general case.
\begin{coro}\label{gohcurve} Assume $\mathbb{G}$ is a Carnot group of step 3, with a Carnot algebra
$\flat=V^1\oplus V^2\oplus V^3$. Let $W$ be a lower-dimensional
subspace of $V^1$.

{\rm(1)} If $[W,V^2]\subsetneq V^3$, then any curve
$(H^1(0)\ni)\gamma^1\subset W$ is a Goh curve.

{\rm(2)} If $\flat$ is a free Carnot algebra, then any curve
$(H^1(0)\ni)\gamma^1\subset W$ is a Goh curve.
\end{coro}
\proof (1) By assumption $\dim (V^3\backslash[W,V^2])\geq1$. For any
curve $(H^1(0)\ni)\gamma^1\subset W$ choosing
$\lambda^1=0,\lambda^2=0$ and $\lambda^3\neq0$ annihilating
$[W,V^2]$, we conclude that $(\gamma^1,\tilde{\lambda})$ is a Goh
extremal. The statement of (2) follows from (1) and the fact that
for a free Carnot algebra, $[W,V^2]\subsetneq V^3$ always holds when
$\dim W<\dim V^1$.
\endproof
Corollary \ref{gohcurve} implies that each Carnot group of step 3
admits Goh curves. The following necessary condition is not used in
this paper, but we include it for completeness.
\begin{prop}\label{GLJ}
Let $\gamma^1, \tilde{\lambda}=(\lambda^1,\lambda^2,\lambda^3)$ be
as in Proposition \ref{goh}. Then
\begin{equation}\label{eq:GLJ}
\lambda^3\left[a,[a,\dot{\gamma}^1(t)]\right]\geq0\textrm{ for a.e.
}t\in[0,1] \textrm{ and any }a\in V^1
\end{equation}
and
\begin{equation}\label{eq:GLJ1}
 \tilde{\lambda}{\rm
d}^2_{\gamma^1}\mathscr{E}(\phi,\phi)=\int^1_0\lambda^3\left[\phi,[\phi,\dot{\gamma}^1]\right]dt\geq0,\quad
\phi\in{\ker\left({\rm d}_{\gamma^1}\mathscr{E}\right)},
\end{equation}
changing $\tilde{\lambda}$ to $-\tilde{\lambda}$ if necessary.
\end{prop}
\proof It suffices to  prove that (\ref{eq:GLJ}) holds for all
Lebesgue points of $\dot{\gamma}^1$.

Let $\bar{\tau}\in[0,1]$ be a Lebesgue point of $\dot{\gamma}^1$.
Assume that
$\lambda^3\left[\bar{a},[\bar{a},\dot{\gamma}^1(\bar{\tau})\right]<0$
for some $\bar{a}\in V^1$. We choose a suitable basis of $V^1$ to
diagonalize the quadratic form
$$
\lambda^3\left[a,[a,\dot{\gamma}^1(\bar{\tau})\right]=\sum_{j=1}^{n_1}\sigma_i(x^j)^2,\quad
a=(x^1,\cdots,x^{n_1})
$$
with at least one term $\sigma_i<0$. For any smooth
$x:\mathbb{R}\rightarrow \mathbb{R}^{n_1}$ with support in [0,1],
let
$$
\phi_x(t)=\left(\underbrace{0,\cdots,0}_{(i-1) \textrm{
terms}},x(t),0,\cdots,0\right),
$$
then we have
\begin{equation}\label{eq:GLJ2}
\tilde{\lambda}{\rm
d}^2_{\gamma^1}\mathscr{E}(\phi_x,\phi_x)=\sigma_i\int_0^1x^2(t)dt<0.
\end{equation}
Denote by $\Pi$ the set of all smooth mappings $x:[0,1]\rightarrow
\mathbb{R}^{n_1}$ with support in [0,1] and satisfying
(\ref{eq:GLJ2}). $\Pi$ is infinitely dimensional, so is
$\left\{\phi_x:\phi_x\in\ker\left({\rm
d}_{\gamma^1}\mathscr{E}\right),x\in\Pi\right\}$. It is a
contraction by Proposition \ref{quadratic1}.
\endproof
(\ref{eq:GLJ}) and (\ref{eq:GLJ1}) are called Generalized
Legendre-Jacobi condition.

\section{Subriemannian Geodesics of step 3}
In this section we assume $\mathbb{G}$ is of step$=3.$
\begin{lemma}\label{line}
Any line (or its interval) through 0 is a normal geodesic.
\end{lemma}
\proof For $\gamma^1(t)=Ctv_0,$ where $C$ a constant and
$v_0\in{V^1}$, we take $\lambda^2=0,\lambda^3=0$, then
(\ref{normal1}) holds for any $\phi\in H^1(0)$ with $\phi(1)=0$.
\endproof

The following result on rank 2 case is well known, see e.g.
\cite{ABCK}. For completeness we give a self-contained proof.

\begin{theorem}[Rank 2 Case]\label{rank2} Let $G$ be of rank 2, i.e., $n_1=2$. Assume $V^1={\rm span}\{e_1,e_2\}$.

{\rm(1)} In the case of the Engel group whose algebra $\flat={\rm
span}\{e_1,e_2\}\oplus{\rm span}\{e_3\}\oplus{\rm span}\{e_4\}$ with
$[e_1,e_2]=e_3, [e_1,e_3]=e_4$, there is a unique arc-length
parameterized singular geodesic $\gamma^1$ which is normal and tangent
to $e_2$, that is, $\gamma^1(t)=te_2$.

{\rm (2)} To the other case where $\flat={\rm
span}\{e_1,e_2\}\oplus{\rm span}\{e_3\}\oplus{\rm span}\{e_4,e_5\}$
with $[e_1,e_2]=e_3, [e_1,e_3]=e_4, [e_2,e_3]=e_5$, singular geodesics are exactly those lines (or their intervals) in $V^1$ through the origin.
\end{theorem}

\proof Let $\gamma^1$ be a singular geodesic parameterized proportionally to arc-length and satisfying
\eqref{singular} with some $\tilde{\lambda}=(\lambda^1,\lambda^2,\lambda^3)$. From \eqref{eq:goh} $\lambda^2=0$. This together with \eqref{singular} implies
that $\lambda^3\neq0$.

We claim that $\lambda^1$ must be in a line. In fact, if this is not true, there must exist $t_1, t_2\in(0, 1]$ such that $V^1=\textrm{span}\{\gamma^1(t_1),\gamma^1(t_2)\}$ which together with \eqref{eq:goh} implies that $\lambda^3=0$
because $[V^1, V^2] = V^3$. A contradiction! So $\gamma^1(t)=(c^1e_1+c^2e_2)t$ for constants $c^1, c^2$ with $(c^1)^2+(c^2)^2\neq0$. In Proposition \ref{lift} we see they are shortest subriemannian geodesics.

(1) If $\flat$ is the Engel algebra, from $\lambda^3\neq0$ and $\lambda^3[c^1e_1t + c^2e_2t, e_3]=0,\forall t\in[0,1]$ we get $c^1=0$, since $[e_1, e_3] = e_4$ and $[e_2, e_3]=0$. From \eqref{F2} and \eqref{F3} by direct computation
we verify that $\gamma^1(t)=c^2e_2t$ $(c^2\neq0)$ is a singular curve (choosing e.g. $\tilde{\lambda}= (0, 0, 1)$).  By Lemma \ref{line} it is normal.

(2) When $\flat$ is the free case, by direct computation we
have $\textrm{Im}\left({\rm d}_{\gamma^1}\mathscr{E}\right)\neq\flat$. In fact, by \eqref{F3} for any $\phi\in H^1(0)$ there exists a constant $\delta$ such that ${\rm d}_{\gamma^1}F^3(\phi)=\delta(c^1e_4+c^2e_5)$. So all $\gamma^1(t)=(c^1e_1+c^2e_2)t$ are singular geodesics. By Lemma \ref{line} they are also normal. \endproof

\begin{lemma}\label{singulartonormal}
Let $\gamma^1$ (or its horizontal lift) be a subriemannian geodesic
in $\mathbb{G}$ and be contained in a lower-dimensional subspace
$W\subset V^1$. If $\gamma^1$ is a normal geodesic in the Carnot
subgroup $\bar{\mathbb{G}}(W)$ of step 2 or 3, then $\gamma^1$ is
also normal in $\mathbb{G}$.
\end{lemma}
\proof Assume $\bar{\mathbb{G}}(W)$ has step 3. Because
$\gamma^1\subset W$ is normal in $\bar{\mathbb{G}}(W)$, by
(\ref{normal1}) there exist $\mu\in{[W,W]^\star}$ and
$\nu\in[W,[W,W]]^\star$ such that
$$
\int_0^1<\dot{\gamma}^1,\dot{\phi}>+\mu[\phi,\dot{\gamma}]+\nu\left[\gamma^1
-\frac{1}{2}\gamma^1(1),\left[\phi,\dot{\gamma}^1\right]\right]dt=0
$$
holds for any $\phi\in H^1([0,1],W)$ with $\phi(0)=\phi(1)=0$. Let
$V^2$ (resp. $V^3$) be orthogonally decomposed as $[W,W]\oplus U^2$
(resp. $[W,[W,W]]\oplus U^3$). Now we take
$\lambda^2=\mu\in{\left(V^2\right)^\star},\lambda^3=\nu\in{\left(V^3\right)^\star}$,
that is, we extend $\mu$ (resp. $\nu$) to $V^2$ (resp. $V^3$) by
annihilating $U^2$ (resp. $U^3$). It is obvious that
$(\lambda^2,\lambda^3)$ satisfies (\ref{normal1}) for any $\phi\in
H^1([0,1],V^1)$ with $\phi(0)=\phi(1)=0$. The case when
$\bar{\mathbb{G}}(W)$ has step 2 is similar.
\endproof

\begin{theorem}[General Case]\label{general} All subriemannian minimizers in $\mathbb{G}$ are
normal.
\end{theorem}
\proof (1) Let $\gamma^1$ be a singular geodesic which is
parameterized proportionally to arc-length and satisfies
(\ref{singular}) for some
$\tilde{\lambda}=(\lambda^1,\lambda^2,\lambda^3)$.

From (\ref{eq:goh}) and (\ref{singular}) we have $\lambda^3\neq0$
because $[V^1,V^1]=V^2,[V^1,V^2]=V^3$. We claim that $\gamma^1$ is
contained in a lower-dimensional subspace $W$ of $V^1$. Otherwise,
there are $t_1,\cdots,t_{n_1}\in{(0,1]}$ such that
$V^1=\textrm{span}\{\gamma^1(t_1),\cdots,\gamma^1(t_{n_1})\}$ which
together with (\ref{eq:goh}) implies $\lambda^3=0$.

Thus $\gamma_1$ (or its horizontal lift) is a subriemannian geodesic
in the Carnot subgroup $\bar{\mathbb{G}}(W)$ whose algebra is
$\bar{\flat}=W\oplus[W,W]\oplus[W,[W,W]]$ or $W\oplus[W,W]$ or $W$.

(2) If $\bar{\flat}=W$, this implies that the horizontal lift of
$\gamma^1$ is itself. By Proposition \ref{lift} the line through 0
and $\gamma^1(t_0)$ for any $t_0\in{(0,1]}$ is the shortest
subriemannian geodesic. So $\gamma^1$ must be an interval of a line
through 0.

(3) If $\bar{\flat}=W\oplus[W,W]$, then $\gamma^1$ is a
subriemannian geodesic in a Carnot group of step 2. So $\gamma^1$ is
normal in $\bar{\mathbb{G}}(W)$. From Lemma \ref{singulartonormal},
$\gamma^1$ is also normal in $\mathbb{G}$.

(4) If $\bar{\flat}=W\oplus[W,W]\oplus[W,[W,W]]$, then  $\gamma^1$
is also a subriemannian minimizer in $\bar{\mathbb{G}}(W)$. If
$\gamma^1$ is regular in $\bar{\mathbb{G}}(W)$, then by Lemma
\ref{singulartonormal} $\gamma^1$ is also normal in $\mathbb{G}$. If
$\gamma^1$ is a singular geodesic in $\bar{\mathbb{G}}(W)$ and $\dim
W\geq3$, we repeat the procedure from step (1), with $\mathbb{G}$
(resp. $\flat)$ replaced by $\bar{\mathbb{G}}(W)$ (resp.
$\bar{\flat})$.

By finite steps we arrive at the case of rank 2. Our statement
follows from Theorem \ref{rank2}.
\endproof
\begin{remark} The smoothness of subriemannian geodesics is very close to the regularity of the subriemannian
distance. In fact, the pointwise smoothness of the subriemannian
distance depends on the strict normalness and uniqueness of
subriemannian geodesics. The subanalyticity of the subriemannian
distance (or sphere) was usually derived from the exclusivity of Goh
curves. We refer the readers to \cite{Ja,AG} and references therein
for this topic. In our case of step 3, as pointed out in Corollary
\ref{gohcurve}, there typically exist Goh curves which are smooth
even normal if they are shortest. Theorem 10 in \cite{AG} proved
that the subriemannian distances of free Carnot groups of step 3 are
not subanalytic.
\end{remark}

\textit{Acknowledgments. Part of this work was done when the first
author visited Department of Mathematics, University of Notre Dame.
He would thank Professor Jianguo Cao for his help and thank the
staff for their hospitality. We also thank the referee for useful
comments and suggestions.}
\bibliographystyle{amsplain}

\end{document}